\documentclass[11pt,francais]{article}
\usepackage{amsmath}\usepackage{epsf,amsfonts,amsthm}\usepackage{amscd,mathptmx,amssymb}
\usepackage{xcolor,epic,eepic}\usepackage{epsfig}
\usepackage{fontenc,indentfirst, delarray,amsfonts,amsmath,amssymb}
\usepackage{rotating}
\usepackage[T1]{fontenc}

\newcommand{\bee}{\begin{enumerate}}
\newcommand{\eee}{\end{enumerate}}
\newcommand{\benn}{\begin{equation*}}
\newcommand{\eenn}{\end{equation*}}
\newcommand{\be}{\begin{equation}}
\newcommand{\ee}{\end{equation}}
\newcommand{\bean}{\begin{eqnarray}}
\newcommand{\eean}{\end{eqnarray}}
\newcommand{\bea}{\begin{eqnarray*}}
\newcommand{\eea}{\end{eqnarray*}}
\newcommand{\w}{\wedge}

\newcommand{\Ci}{C^{\infty}}
\newcommand{\E}{\ell}

\newcommand{\R}{\mathbb{R}}

\newcommand{\n}{\nabla}

\newcommand{\lp}{\left(}
\newcommand{\rp}{\right)}

\newcommand{\op}[1]{\!\!\mathop{\rm ~#1}\nolimits}

\newcommand{\bqa}{\begin{eqnarray}}
\newcommand\eqa {\end{eqnarray}}
\newcommand{\beq}{\begin{eqnarray}}
\newcommand{\beqn}{\begin{eqnarray}\nonumber}
\newcommand{\eeq}{\end{eqnarray}}

\newcommand{\wcK}{{\widetilde{\cal K}}}
\newcommand{\wA}{{\widetilde{A}}}
\newcommand{\wzp}{{\widetilde{\zp}}}
\newcommand{\wa}{{\widetilde{a}}}
\newcommand{\wcA}{{\widetilde{\cal A}}}
\newcommand{\cZ}{{\cal Z}}
\newcommand{\sZ}{{\textsf{Z}}}

 \newcommand{\cK}{{\cal K}}

 \newcommand{\cA}{{\cal A}}

 \newcommand{\cI}{{\cal I}}

\def\xd{\mathrm{d}}

   \newcommand{\zf}{{\phi}}
   \newcommand{\zvf}{{\varphi}}

\mathchardef\za="710B  
\mathchardef\zb="710C  
\mathchardef\zg="710D  
\mathchardef\zd="710E  
\mathchardef\zve="710F 
\mathchardef\zz="7110  
\mathchardef\zh="7111  
\mathchardef\zy="7112 

\mathchardef\zi="7113  
\mathchardef\zk="7114  
\mathchardef\zl="7115  
\mathchardef\zm="7116  
\mathchardef\zn="7117  
\mathchardef\zx="7118  
\mathchardef\zp="7119  
\mathchardef\zr="711A  
\mathchardef\zs="711B  
\mathchardef\zt="711C  
\mathchardef\zu="711D  
\mathchardef\zf="711E 
\mathchardef\zq="711F  
\mathchardef\zc="7120  
\mathchardef\zw="7121  
\mathchardef\ze="7122  
\mathchardef\zvy="7123  
\mathchardef\zvw="7124  
\mathchardef\zvr="7125 
\mathchardef\zvs="7126 
\mathchardef\zvf="7127  
\mathchardef\zG="7000  
\mathchardef\zD="7001  
\mathchardef\zY="7002  
\mathchardef\zL="7003  
\mathchardef\zX="7004  
\mathchardef\zP="7005  
\mathchardef\zS="7006  
\mathchardef\zU="7007  
\mathchardef\zF="7008  
\mathchardef\zW="700A  

\newcommand{\cyclic}{\mathop{\kern0.9ex{{+}
 \kern-2.15ex\raise-.25ex\hbox{\Large\hbox{$\circlearrowright$}}}}\limits}
\newtheorem{theo}{Theorem}
\newtheorem{thm}{Theorem}
\newtheorem{prop}{Proposition}
\newtheorem{lem}{Lemma}

\theoremstyle{definition}

\newtheorem{rem}[thm]{Remark}

\begin{document}
\title{Geometric structures encoded in the\\ Lie structure of an Atiyah algebroid}
\author{Janusz Grabowski\footnote{The research of J.~Grabowski was supported by the Polish
Ministry of Science and Higher Education under the grant No. N201 005 31/0115.}, Alexei Kotov, Norbert Poncin
\footnote{The research of N. Poncin was supported by UL-grant SGQnp2008.}
\date{7 May 2009}}
\maketitle

\begin{abstract} We investigate Atiyah algebroids, i.e. the infinitesimal objects
of principal bundles, from the viewpoint of Lie algebraic approach
to space. First we show that if the Lie algebras of smooth
sections of two Atiyah algebroids are isomorphic, then the
corresponding base manifolds are necessarily diffeomorphic.
Further, we give two characterizations of the isomorphisms of the
Lie algebras of sections for Atiyah algebroids associated to
principle bundles with semisimple structure groups. For instance
we prove that in the semisimple case the Lie algebras of sections
are isomorphic if and only if the corresponding Lie algebroids
are, or, as well, if and only if the integrating principal bundles
are locally diffeomorphic. Finally, we apply these results to
describe the isomorphisms of sections in the case of reductive
structure groups -- surprisingly enough they are no longer
determined by vector bundle isomorphisms and involve divergences
on the base manifolds.\end{abstract}

\vspace{2mm} \noindent {\bf MSC 2000}: 17B40, 55R10, 57R50, 58H05 (primary), 16S60, 17B20, 53D17
(secondary)\medskip


\section{Introduction}

The concept of groupoid was introduced in 1926 by the German mathematician Heinrich Brandt. A groupoid is a
small category in which any morphism is invertible. Topological and differential groupoids go back to Charles
Ehresmann \cite{EhresmannLieGroupoids} and are transitive Lie groupoids in the sense of A. Kumpera and A.
Weinstein: Lie groupoids are groupoids whose classes $\zG_0$ of objects and $\zG_1$ of morphisms are not only
sets, but manifolds, the source and target maps $s$ and $t$ are submersions, and all operations are smooth;
such a groupoid is called transitive, if any two objects are related by a morphism, i.e. if
$(s,t):\zG_1\to\zG_0\times\zG_0$ is surjective. The gauge groupoid of a principal $G$-bundle $P$, $\zG_0=P/G$
and $\zG_1=(P\times P)/G$, where the $G$-action on pairs is componentwise, is the prototype of a transitive
Lie groupoid; actually, any transitive Lie groupoid can be viewed as the gauge groupoid of a principal bundle
\cite{LibermannTransLad}.

The first-order invariants of principal bundles or transitive Lie
groupoids $P(M,\zp,G)$ are {\it Atiyah sequences} of vector
bundles
\[ 0\to
K:=P\times_G{\frak g}\to A\stackrel{\zp_*}{\to}TM\to 0,\] together
with the corresponding sequences of modules of sections
\[0\to {\cal K}\to {\cal A}{\to}{\frak X}(M)\to 0\,.
\]
They were introduced by M. F. Atiyah \cite{AtiyahSequence} in order to investigate the existence of complex
analytic connections in fiber bundles. Atiyah sequences were referred to as {\it Atiyah algebroids} after the
introduction of Lie algebroids by J. Pradines \cite{PradinesLad} in order to grasp the structure of the
infinitesimal objects that correspond to Lie groupoids; Lie algebroids unify several well-known passages to
the infinitesimal level: from foliation to distribution, from Lie group action to Lie algebra action, from
principle bundle to Atiyah sequence, from symplectic Lie groupoid to Poisson manifold, etc. Atiyah sequences
or algebroids are in particular transitive Lie algebroids, in the sense that their anchor map is surjective.

Let us also mention that Atiyah algebroids naturally appear in
Theoretical Mechanics. Indeed, subsequently to A. Weinstein's work
on the unification of internal and external symmetry
\cite{WeinsteinInternalExternalSymm} and on the groupoid approach
to Mechanics \cite{We}, Lagrangian functions on Lie algebroids
were investigated; if we consider a Lagrangian with symmetries on
a configuration space that is a principal $G$-bundle $P$, i.e. a
Lagrangian that is invariant under the action of the structure Lie
group $G$, then the Lagrangian function is defined on $A:=TP/G$,
i.e. on the Atiyah algebroid associated with this principal bundle
\cite{LéonMarreroMartinezDynamicsLieAlg,
\cite{MestdagAtiyahAlgMechanics}.} This gives rise to different
Lie algebroid generalizations of the Lagrange and Hamilton
formalisms \cite{Mar}, \cite{GGU}, \cite{GG}.\medskip

In this work, {\it we investigate Atiyah algebroids from the standpoint of Lie algebraic approach to space}.
The general result \cite[Theorem 8]{GrabowskiGrabowska} implies that if the Lie algebras of smooth sections of
two Atiyah algebroids $(A_i,[-,-]_i,\zp_{i*})$, $i\in\{1,2\}$, over two differentiable manifolds $M_i$ are
isomorphic, then the base manifolds $M_i$ are necessarily diffeomorphic, Section 3, Theorem 1; see also
\cite{GrabPoncinDiffOpLineBundles}. We go further and characterize the isomorphisms of the Lie algebras of
smooth sections $({\cal A}_i,[-,-]_i)$, $i\in\{1,2\}$, for Atiyah algebroids associated with principal bundles
with semisimple structure groups, Section 4, Theorem 2. {\it They turn out to be associated with Lie algebroid
isomorphisms.} To obtain these upshots, we identify the maximal finite-codimensional Lie ideals of the kernel
${\cal K}$ of the corresponding Atiyah sequence of modules and Lie algebras, Section 4, Theorem 3, and
describe the elements of ${\cal K}$ and ${\cal A}$, which vanish at a given point, in pure Lie algebraic
terms, Section 4, Theorems 5 and 6. Next, we prove that Lie algebra isomorphisms between sections come from
vector bundle isomorphisms and characterize them. Note that the assumption of semisimplicity is essential to
prove that a Lie algebra isomorphism for sections is implemented by an isomorphism of the vector bundles, as
shows the case of the Atiyah algebroid $TM\times\R$ of first-order differential operators
\cite{GrabPoncinAuto}. Combining the semisimple case with the case of first-order differential operators, we
describe the isomorphisms for reductive structure groups. They no longer come from vector bundle isomorphisms,
as first-order components associated with divergences appear in the picture. Denoting by $\sZ$ the subbundle
$\sZ\subset K$ associated with the center $Z{\frak g}$ of the Lie algebra $\frak g$ of the structure group
$G$, our main result can be stated as follows.

\begin{theo}\label{t0}Let \ $\zF:\cA_1\rightarrow\cA_2$ be an isomorphism of the Lie algebras $(\cA_i,[-,-]_i)$ of Atiyah
algebroids $$0\to K_i\to A_i\stackrel{\zp_{i*}}{\to} TM_i\to 0$$
with connected reductive structure groups $G_i$ over connected
manifolds $M_i$, $i\in\{1,2\}$. Then $\zF$ is a unique composition
$\zF=\zF_0\circ\zF_1$, where $\zF_0$ is an isomorphism of the Lie
algebras $(\cA_i,[-,-]_i)$ associated with a Lie algebroid
isomorphism and $\zF_1$ is an automorphism of the Lie algebra
$(\cA_1,[-,-]_1)$ of the form
\be\label{11}\zF_1(a)=a+\op{div}(\zp_{1*}a)\cdot{r}\,, \ee
where $\op{div}:{\frak X}(M_1)\rightarrow C^\infty(M_1)$ is a divergence operator on $M_1$ and ${r}$ is a
section of $\sZ_1$ represented by the fundamental vector field of an element $r\in Z{\frak g}_1$.
\end{theo}
Let us mention that {\it algebraic characterization of space} can
be traced back to Gel'fand and Kolmogoroff -- description of
isomorphisms between associative algebras of continuous functions
on compact sets -- and is concerned with characterization of
diverse geometric structures by means of various associative or
Lie algebras growing on them. It is well known that isomorphisms
of the associative algebras of smooth functions living on second
countable manifolds are implemented by diffeomorphisms of the
underlying manifolds (for the general case, including manifolds
which are not second countable nor paracompact, see \cite{Gra05,
Mrc05}). The classical result by Pursell and Shanks \cite{PS},
which states that the Lie algebra structure of the space of
compactly supported vector fields characterizes the differential
structure of the underlying manifold, is the starting point of a
multitude of papers: Koriyama, Maeda, Omori (other complete and
transitive Lie algebras of vector fields), see e.g. \cite{Omori},
Amemiya, Masuda, Shiga, Duistermaat, Singer, Grabowski, Poncin
(differential and pseudodifferential operators)
\cite{AMS,DuistermaatSinger,GrabPoncinAuto}, Abe, Atkin,
Grabowski, Fukui, Tomita, Hauser, Müller, Rybicki (special
geometric situations), Amemiya, Grabowski (real and analytic
cases), see e.g. \cite{GrabowskiIsomIdeal}, Skryabin (modular Lie
algebras) \cite{Skryabin}, and Grabowska, Grabowski (Lie algebras
associated with Lie algebroids) \cite{GrabowskiGrabowska}.

\section{Atiyah algebroids, Lie algebra bundles, and Lie algebroids}\label{RemarksAlds}

To ensure independent readability of the present text, we recall some facts about Atiyah and Lie algebroids.

\begin{rem} In this article all manifolds are smooth second countable
paracompact Hausdorff manifolds of finite and nonzero dimension.\end{rem}

\subsection{Atiyah algebroids}

Let $P(M,\zp,G)$ be a principal $G$-bundle $\zp:P\to M$. Set ${\frak g}=\op{Lie}(G)$ and denote by $r_g:P\ni
u\to u.g\in P$ the diffeomorphism that is induced by the right action of $g\in G$ on $P$. Over any point $m\in
M$, we define an equivalence relation of tangent vectors to $P$. If $\zp(u)=m$, $g\in G$, $X_u\in T_uP$, and
$X_{u.g}\in T_{u.g}P$, the vectors $X_u$ and $X_{u.g}$ are said to be equivalent if and only if
$X_{u.g}=(T_ur_{g})(X_u)$, where $T_ur_g$ is the tangent isomorphism. The equivalence classes of this relation
form a vector space $A_m\simeq T_uP$, and the disjoint union $A=\bigsqcup_{m\in M}A_m$ is a {\it vector
bundle} that is---as quite easily seen---locally diffeomorphic to $A\vert_U\simeq TU\times {\frak g}$, where
$U$ is a chart domain of $M$ over which $P$ is trivial. As aforementioned, the vector bundle $A$ is often
denoted by $TP/G$.\medskip

We use the same notations as above. Since $\zp\circ r_g=\zp$, it is clear that the image of a vector $X_u$ by
the surjection $T_u\zp:T_uP\to T_mM$ does not depend on the representative $X_u$ of the class $[X_u]\in A_m$,
but only on the class itself. Hence, we get a well-defined surjection $\zp_{*m}:A_m\to T_mM$, as well as a
surjective bundle map $\zp_*:A\to TM$ over the identity. This map will be the {\it anchor} of the Atiyah
algebroid $(A,[-,-],\zp_*)$ associated with the principal bundle $P(M,\zp,G)$.\medskip

To define the Lie bracket $[-,-]$ on the space ${\cal A}:=\zG(A)$ of smooth sections of $A$, consider the
short exact Atiyah sequence of vector bundles and bundle maps
\be 0\to K=P\times_G{\frak g}\to
A\stackrel{\zp_*}{\to}TM\to 0.\label{AtiyahSequence}
\ee
Let us roughly explain why the kernel
$K:=\op{ker}\zp_*$ coincides with the associated vector bundle $P\times_G{\frak g}$ over $M$. Since at each
point $u\in P$ the vertical tangent vectors coincide with the fundamental vectors $X_u^h$, $h\in{\frak g}$, of
the action, we have
$$K=\op{ker}\zp_*=\bigsqcup_{m\in M}\{[X_u]:T_u\zp\;
X_u=0,\zp(u)=m\}=\bigsqcup_{m\in M}\{[X_u^h]:\zp(u)=m,h\in{\frak g}\}.$$ As for each $u\in P$, the map
$h\in{\frak g}\to X_u^h\in V_u$ is a vector space isomorphism between the space ${\frak g}$ and the space
$V_u$ of vertical vectors at $u$, and as $$T_ur_g\; X_u^h=X_{u\cdot g}^{\op{Ad}(g^{-1})h},$$ where $\op{Ad}$
is the adjoint representation of $G$ on the vector space ${\frak g}$, the afore-detailed equivalence relation
identifies $(u,h)\simeq (u.g,\op{Ad}(g^{-1})h)$. Hence, the kernel $K$ is actually the associated vector
bundle $P\times_G{\frak g}$.

As soft sheaves over paracompact spaces are acyclic, a short exact
sequence of vector bundles over $M$ induces a short exact sequence
of the $\Ci(M)$-modules of sections. Since the module ${\frak
X}_G(P)$ of $G$-invariant vector fields of $P$ visibly coincides
with the module ${\cal A}$ of sections of $A$, this new sequence
is \be 0\to {\cal K}:=\zG(K)\simeq\Ci(P,{\frak g})^G\to {\cal
A}:=\zG(A)\simeq {\frak X}_G(P)\stackrel{\zp_*}{\to}{\cal
V}:={\frak X}(M)\to 0,\label{AtiyahSequenceSections} \ee where
$\Ci(P,{\frak g})^G$ is the module of $G$-equivariant smooth
functions from $P$ to ${\frak g}$, and where ${\frak X}(M)$ is the
module of vector fields of $M$. Sequence
(\ref{AtiyahSequenceSections}) is also a short exact sequence of
Lie algebras. The {\it Lie bracket} $[-,-]$ of ${\cal A}$ is of
course implemented by the Lie algebra structure of ${\frak
X}_G(P)$. Indeed, this subspace is a Lie subalgebra of ${\frak
X}(P)$, since a $G$-invariant vector field of $P$ is a vector
field that is $r_g$-related to itself, for all $g\in G$. The
compatibility property of the bracket of vector fields of $P$ with
the $\Ci(P)$-module structure, entails the corresponding property
of the bracket of ${\cal A}$ with the $\Ci(M)$-module structure.
As the anchor $\zp_*$ of the {\it Atiyah algebroid}
$(A,[-,-],\zp_*)$ is automatically a Lie algebra homomorphism, the
kernel ${\cal K}$ of $\zp_*$ is a Lie subalgebra of ${\cal A}$ and
even a Lie ideal. The {Atiyah algebroid} $(A,[-,-],\zp_*)$ is
therefore an example of a {\it transitive Lie algebroid} and the
exact sequence (\ref{AtiyahSequence}) is an exact sequence of
morphisms of Lie algebroids.

\medskip
It is well-known and easily checked that the Lie bracket of any Lie algebroid is local and that a Lie
algebroid thus restricts to a Lie algebroid over any open subset of the base of the initial bundle. Let now
$U\subset M$ be simultaneously a chart domain of $M$ and a trivialization domain of $P$\,: $A\vert_U\simeq
TU\times{\frak g}$ and $\zG(A\vert_U)\simeq {\frak X}(U)\times\Ci(U,{\frak g})$. As $\zp\vert_U:P\vert_U\simeq
U\times G\ni (m,g)\to m\in U$, one immediately sees that $\zp_*\vert_U:TU\times{\frak g}\to TU$ is the
projection onto the first factor. Hence, if $(v_1,\zg_1),(v_2,\zg_2)\in{\frak X}(U)\times\Ci(U,{\frak g})$ are
sections of $A$ over $U$, the Lie algebra homomorphism property of the restriction algebroid shows that the
first component of their bracket is $[v_1,v_2]_{{\frak X(U)}}$. Further, it follows from the compatibility
condition between the module and the Lie structures in a Lie algebroid, that in the kernel ${\cal K}$ the Lie
bracket is $\Ci(M)$-bilinear, and thus provides a bracket in each fiber $K_m\simeq{\frak g}$ of
$K=P\times_G{\frak g}$. This bracket $[k_m,k'_m]=[k,k']_m$, $k,k'\in{\cal K}$, is actually the bracket of
${\frak g}$ \cite{KMacSecondBook}, so that $K$ becomes a Lie algebra bundle ({\small LAB}). Hence, the value
at $m\in U$ of the bracket $[(0,\zg_1),(0,\zg_2)]$ of the elements $(0,\zg_i)\in\zG(K\vert_U)$ coincides with
the bracket $[\zg_{1m},\zg_{2m}]_{\frak g}$. Since the module-Lie compatibility condition shows that for any
$f\in\Ci(U)$, we have
$$\Ci(U,{\frak g})\ni
[(v_1,0),(0,f\zg_2)]=f[(v_1,0),(0,\zg_2)]+v_1(f)\zg_2,$$ we eventually realize that
\be\label{AtiyahAldBracketLocal} [(v_1,\zg_1),(v_2,\zg_2)] =([v_1,v_2]_{{\frak X}(U)},[\zg_1,\zg_2]_{\frak
g}+v_1(\zg_2)-v_2(\zg_1)),\ee where the ${\frak g}$-bracket is computed pointwise, although the exact proof of
the last equation is highly nontrivial. This semidirect product and the aforementioned projection onto the
first factor define on $TU\times{\frak g}$ a Lie algebroid structure that is called {\it trivial Lie
algebroid} on $U$ with structure algebra ${\frak g}$ (the direct product would not define a Lie algebroid).

\subsection{Lie algebra bundles}

In the following, $K$ denotes an arbitrary {\small LAB} with
typical fiber $\frak g$ over a manifold $M$. First remember, see
\cite[Prop. 3.3.9]{KMacSecondBook}, that if $\frak h$ denotes a
Lie subalgebra of $\frak g$, there is a sub-{\small LAB} $H$ of
$K$ with typical fiber $\frak h$, whose {\small LAB}-atlas is
obtained by restriction of the {\small LAB}-atlas of $K$, on the
condition that $\frak h$ is preserved by any automorphism of
$\frak g$ -- so that the restriction of the transition cocycle of
$K$ provides a transition cocycle of $H$. If $\frak h$ is the
center $Z\frak g$ of $\frak g$ (resp. the derived ideal $[\frak
g,\frak g]$ of $\frak g$), the corresponding sub-{\small LAB} $ZK$
(resp. $[K,K]$) is called the {\it center sub-{\small LAB}} of $K$
(resp. the {\it derived sub-{\small LAB}} of $K$) and it is
denoted shortly by $\sZ$ (resp. $K^{(1)}$).\medskip

It is known and easily seen that
\be\zG(\sZ)=Z\,\zG(K)=:{\cZ}\label{SectionsCenter}\ee ( resp. that
\be\label{SectionsDerived}\zG([K,K])=[\zG(K),\zG(K)]=[{\cal
K},{\cal K}]\;\;).\ee To understand for instance why
$\zG([K,K])\subset [{\cal K},{\cal K}]$, denote by $g_i$ some
global generators of the module $\zG(K)$ over $\Ci(M)$, observe
that for any open subset $U\subset M$ any element of
$[\zG(K\vert_U),\zG(K\vert_U)]$ reads \be\sum
[k'_U,k''_U]=\sum[k'_U,\sum_ik_U^{''\,i}\,g_i\vert_U]=\sum_i[\sum
k_U^{''\,i}k'_U,g_i\vert_U]=:\sum_i[k_U^{i},g_i\vert_U],\label{BrackSectKLocal1}\ee
with $k_U^{i}\in\zG(K\vert_U)$. Take now a partition of unity
$(U_{\za},\zvf_{\za})$ that is subordinated to a cover by local
trivializations of $K$ and let $\zd\in\zG([K,K])$. Since
$[K,K]\vert_{U_{\za}}$ is diffeomorphic to $U_{\za}\times[\frak
g,\frak g]$, where $\frak g$ is finite-dimensional, we have
\be\zG([K,K]\vert_{U_{\za}})=\Ci(U_{\za},[\frak g,\frak
g])=[\Ci(U_{\za},\frak g),\Ci(U_{\za},\frak
g)]=[\zG(K\vert_{U_{\za}}),\zG(K\vert_{U_{\za}})].\label{BrackSectKLocal2}\ee
When combining the upshots (\ref{BrackSectKLocal2}) and
(\ref{BrackSectKLocal1}), we get
$$\zd=\sum_{\za}\zvf_{\za}\zd\vert_{U_{\za}}=\sum_i[\sum_{\za}\zvf_{\za}\,k^{i}_{U_{\za}},g_i],$$
with obvious notations. The conclusion $\zG([K,K])\subset[\zG(K),\zG(K)]$ follows. If $\frak g$ is semisimple,
Equation (\ref{BrackSectKLocal2}) shows that locally we have
$\zG(K\vert_{U_{\za}})=[\zG(K\vert_{U_{\za}}),\zG(K\vert_{U_{\za}})]$. When combining this conclusion again
with Equation (\ref{BrackSectKLocal1}), but now for $\zd\in\zG(K)$, we see just as before that
$\zG(K)\subset[\zG(K),\zG(K)]$. In the semisimple case we thus get \be {\cal K}={\cal K}^{(1)}=[{\cal K},{\cal
K}].\label{KernelSemisimple}\ee

Eventually, a {\small LAB} $K$ is said to be {\it reductive}, if all its fibers are, hence, if its typical
fiber $\frak g$ is. In that case, the structure of the {\small LAB} is, exactly as in Lie algebra theory, \be
K=\sZ\oplus[K,K]\label{KernelReducive}.\ee

\subsection{Lie algebroids}\label{Lads}

If $(L,[-,-],\zr)$ and $(L',[-,-]',\zr')$ are two Lie algebroids
over two manifolds $M$ and $M'$ respectively, a Lie algebroid
morphism $F$ should of course in particular be a vector bundle map
$F:L\to L'$ over some smooth map $f:M\to M'$. Since, if $f$ is not
a diffeomorphism, such a bundle morphism does not necessarily
induce a morphism $F:{\cal L}\to{\cal L}'$ between the
corresponding modules of sections ${\cal L}:=\zG(L)$, ${\cal
L'}:=\zG(L')$, the definition of a Lie algebroid morphism is
nontrivial in the general situation. However, if $f:M\to M'$ is a
diffeomorphism, and especially in the base-preserving case where
$f$ is just identity, a {\it Lie algebroid morphism} can be
naturally defined as a vector bundle map $F:L\to L'$ that verifies
$\zr'\circ F=f\circ\zr$ and $F[\E_1,\E_2]=[F\E_1,F\E_2]',$ for any
$\E_1,\E_2\in{\cal L}$. Of course, such a morphism is called a
{\it Lie algebroid isomorphism}, if $F:L\to L'$ is a vector bundle
isomorphism.\medskip

It is well known that any short exact sequence of vector spaces
and linear maps splits (in the infinite-dimensional setting the
result is based upon the axiom of choice). The same is true for a
short exact sequence of vector bundles $$0\to E\stackrel{i}{\to}
F\stackrel{\zr}{\to} G\to 0$$ over a same manifold and vector
bundle maps that cover the identity. In the following, we
systematically assume for simplicity reasons that $E\subset F$ is
a vector subbundle of $F$ and that $i$ is just the inclusion. As
for the mentioned {\it splitting of any short exact sequence of
vector bundles}, it suffices to consider a smooth Riemannian
metric in $F$ and to define the orthogonal vector subspace
$E_m^{\bot}$ of each fiber $E_m\subset F_m$ with respect to that
metric. These orthogonal subspaces then glue smoothly and form a
vector subbundle $E^{\bot}\subset F$, such that $E\oplus
E^{\bot}=F$. It is easily checked that existence of a vector
subbundle of $F$ that is supplementary to $E$ in $F$ is equivalent
to existence of a vector bundle isomorphism $\zt:F\to E\oplus G$,
which nevertheless has to ``work'' according to the conditions
$\zt\circ i=(\op{id}_E, 0)$ and $\zr=\op{pr_2}\circ\zt$, where
$\op{pr_2}$ is the projection onto the second term. A third
equivalent definition of splitting asks for a vector bundle map
$\theta:G\to F$ (resp. $j:F\to E$) that is a right inverse to
$\zr$ (resp. a left inverse to $i$). Let us eventually mention
that a splitting (say $\theta$) of a short exact sequence of
vector bundles over a same base $M$, induces of course in the
natural way a splitting (we will denote it by $\theta$ as well) of
the corresponding short exact sequence of $\Ci(M)$-modules of
sections.

Consider now a transitive Lie algebroid $(L,[-,-],\zr)$ over $M$
and let $$0\to K\stackrel{i}{\to}L\stackrel{\zr}{\to}TM\to 0$$ be
the corresponding short exact sequence of vector bundles -- where
$K$ is known to be a {\small LAB}. It is customary to refer to a
right inverse bundle map or right splitting $\theta$ (resp. left
inverse bundle map or left splitting $j$) as a {\it Lie algebroid
connection} of $L$ (resp. {\it connection reform} of $L$). The
point is of course that if the investigated algebroid is the
Atiyah algebroid $(A,[-,-],\pi_*)$ of a principal bundle
$P(M,\zp,G)$, there is a 1-to-1 correspondence between Lie
algebroid connections of $A$ (resp. connection reforms of $A$) and
connections of the principal bundle $P$ (resp. connection 1-forms
of $P$) in the traditional sense. The {\it curvature $R_{\theta}$
of a transitive Lie algebroid connection $\theta$} measures the
Lie algebra morphism default of $\theta$, i.e. for any vector
fields $X,Y\in{\frak X}(M)$, we set
$$R_{\theta}(X,Y):=[\theta X,\theta Y]-\theta[X,Y]\in {\cal K}.$$
Since $R_{\theta}$ is, as easily seen, $\Ci(M)$-bilinear, it defines a vector bundle map $R_{\theta}:TM\times
TM\to K$.\medskip

An {\it ideal of a transitive Lie algebroid} $$0\to
K\stackrel{i}{\to}L\stackrel{\zr}{\to}TM\to 0$$ is a sub-{\small
LAB} $H$ of $K$, such that, for any $h\in{\cal H}$ and $\E\in{\cal
L}$, we have $[h,\E]\in{\cal H}$. For instance, the {\small LAB}
$K$ itself, its center sub-{\small LAB} $\sZ$ and its derived
sub-{\small LAB} $[K,K]$ are Lie algebroid ideals of $L$.

Let $H$ be any Lie algebroid ideal of $L$ and consider the short
exact sequence of vector bundles \be\label{SesLads}0\to
H\stackrel{\zi}{\to} L\stackrel{p}{\to} L/H\to 0.\ee Since $$0\to
{\cal H}\stackrel{\zi}{\to} {\cal L}\stackrel{p}{\to} \zG(L/H)\to
0$$ is a short exact sequence of $\Ci(M)$-modules, the module
morphism $p$ induces a canonical isomorphism of modules between
${\cal L}/{\cal H}$ and $\zG(L/H)$. This isomorphism allows
transferring the Lie algebra structure $[-,-]^{\;\tilde{}}$ of
${\cal L}/{\cal H}$ to $\zG(L/H)$. As the vector bundle map
$\zr:L\to TM$ factors through the quotient $L/H$,
$\zr=\tilde\zr\circ p$, the triplet
$(L/H,[-,-]^{\;\tilde{}},\tilde \zr)$ is a transitive Lie
algebroid over $M$ -- the {\it quotient Lie algebroid} of $L$ over
the Lie algebroid ideal $H$ -- and $\op{ker}\tilde \zr=K/H$.
Moreover, it is clear that the sequence (\ref{SesLads}) is a short
exact sequence of Lie algebroids and Lie algebroid morphisms.

\section{Lie algebras of Atiyah algebroids}

In the following, we compendiously investigate whether an
isomorphism between the Lie algebras of sections of two Atiyah
algebroids induces a diffeomorphism between the underlying
manifolds.

\medskip
Consider any smooth Lie algebroid $(L,[-,-],\widehat{-}\;)$ over a
smooth manifold $M$, and set, in order to simplify notations,
${\cal N}=\Ci(M)$ and ${\cal L}=\zG(L)$. We say that the algebroid
$L$ is {\it strongly nonsingular}, if ${\cal N}=\widehat{\cal
L}({\cal N}):=\op{span}\{\widehat{\E}(f):\E\in{\cal L},f\in{\cal
N}\}$. As, in view of the Serre-Swan theorem, the ${\cal
N}$-module ${\cal L}$ can be characterized as a projective module
with a finite number of generators, the generalized foliation
spanned by $\widehat{\cal L}$ is finitely generated and the strong
nonsingularity can be characterized in geometric terms as the fact
that there is a finite number of vector fields from $\widehat{\cal
L}$ which do not all vanish at a single point \cite{Gra93}. In
\cite[Theorem 8]{GrabowskiGrabowska} it has been proved that under
this condition isomorphisms between the Lie algebras of Lie
algebroids induce diffeomorphisms of the underlying
manifolds.\medskip

To fix notation for further purposes, let us briefly sketch the proof of this fact for Atiyah algebroids.
Denote by ${I}({\cal N})$ (resp. ${S}({\cal L})$) the set of all maximal finite-codimensional associative
ideals of ${\cal N}$ (resp. the set of all maximal finite-codimensional Lie subalgebras of ${\cal L}$). It was
shown in \cite[Corollary 7]{GrabowskiGrabowska} that, in the strongly nonsingular case, we can pass from the
associative to the Lie setting via the action of vector fields on functions, i.e., more precisely, that the
map $$ {I}({\cal N})\ni J\leftrightarrow {\cal L}_J:=\{\E\in{\cal L}:\widehat{\E}({\cal N})\subset
J\}\in{S}({\cal L})$$ is a bijection. Since the maximal finite-codimensional ideals of ${\cal N}$ ``are''
exactly the points of $M$, i.e. as the map \be M\ni m\leftrightarrow J(m):=\{f\in{\cal N}:f(m)=0\}\in I({\cal
N})\label{Maxfincodimidealsfcts}\ee is bijective \cite[Proposition 3.5]{GrabowskiIsomIdeal}, we get a 1-to-1
correspondence $$ M\ni m\leftrightarrow {\cal L}_m:={\cal L}_{J(m)}=\{\E\in{\cal
L}:\widehat{\E}(f)(m)=0,\forall f\in{\cal N}\}=\{\E\in{\cal L}:\widehat{\E}_m=0\}\in S({\cal L}).$$

Let $(A,[-,-],\zp_*)$ be an Atiyah algebroid over a manifold $M$. We use the above notations; further, we set
for convenience, $\widehat{-}:=\zp_*$. It follows from transitivity that $\widehat{\cal A}={\cal V}:={\frak
X}(M)$, so that $A$ is strongly nonsingular. Since this observation entails that
$$ M\ni m\leftrightarrow {\cal A}_m=\{a\in{\cal
A}:\widehat{a}_m=0\}\in S({\cal A}),$$ the kernel ${\cal K}$ of the Atiyah sequence of $\Ci(M)$-modules and
Lie algebras of sections associated with $A$ reads $${\cal K}=\cap_{m\in M}\{a\in{\cal
A}:\widehat{a}_m=0\}=\cap_{T\in S({\cal A})} T.$$

If $\zF:{\cal A}_1\leftrightarrow {\cal A}_2$ is a Lie algebra isomorphism, we have $\zF({\cal
K}_1)=\cap_{T\in S({\cal A}_1)} \zF(T)={\cal K}_2,$ since $\zF$ maps maximal finite-codimensional Lie
subalgebras into maximal finite-codimensional Lie subalgebras. Hence, $\zF$ induces an isomorphism
$\widetilde{\zF}$ between the quotient Lie algebras ${\cal A}_i/{\cal K}_i\simeq{\cal V}_i$, $i\in\{1,2\}$.
Such an isomorphism however, is implemented by a diffeomorphism $\zvf:M_1\leftrightarrow M_2$ \cite[Corllary
5.8]{GrabowskiIsomIdeal}, $(\widetilde{\zF} v)f=(v(f\circ \zvf))\circ \zvf^{-1}$, for any $v\in{\cal V}_1$ and
any $f\in{\cal N}_2$. Hence, we can formulate the following Pursell-Shanks type theorem which is a particular
case of \cite[Theorem 8]{GrabowskiGrabowska}.

\begin{theo}\label{TheoPursellShanks} Let $(A_i,[-,-]_i,\zp_{i*})$ be a smooth Atiyah algebroid over
a smooth manifold $M_i$, $i\in\{1,2\}$. Any isomorphism $\zF:{\cal A}_1\leftrightarrow {\cal A}_2$ between the
Lie algebras of smooth sections restricts to an isomorphism $\zF:{\cal K}_1\leftrightarrow{\cal K}_2$ between
the kernels of the anchor maps and induces an isomorphism $\widetilde{\zF}:{\cal V}_1\leftrightarrow {\cal
V}_2$ between the Lie algebras of vector fields of the base manifolds, which is implemented by a
diffeomorphism $\zvf:M_1\leftrightarrow M_2$, i.e. $\widetilde{\zF}=\zvf_*$. Furthermore, $\zF({\cal
A}_{1m})={\cal A}_{2\zvf(m)}$.\end{theo}

\begin{proof} We only need prove the last claim. It is clear that $\zF({\cal A}_{1m})={\cal
A}_{2n}$, for some $n\in M_2$. Since the abovementioned isomorphism between ${\cal A}_i/{\cal K}_i$ and ${\cal
V}_i$ is $\widetilde{\zp}_{i*}: {\cal A}_i/{\cal K}_i\ni [a]\leftrightarrow \zp_{i*}a\in{\cal V}_i,$ we have
in fact $\widetilde{\zp}_{2*}\widetilde{\zF}=\zvf_*\widetilde{\zp}_{1*}$, so that $\widehat{\zF
(a)}=\widetilde{\zp}_{2*}\widetilde{\zF}[a]=\zvf_*\widehat{a}.$ When evaluating both sides of this equation at
$\zvf(m)$, we get
$$\widehat{\zF (a)}_{\zvf(m)}=(\zvf_{*m})\widehat{a}_m.$$ Hence the
result. \end{proof}

\section{Isomorphisms of Lie algebras of Atiyah algebroids - semisimple structure
groups}\label{SemisimpleStrucGroups}

In this section we take an interest in a possible characterization of isomorphisms $\zF:{\cal
A}_1\leftrightarrow {\cal A}_2$ of the Lie algebras of smooth sections ${\cal A}_i$ of Atiyah algebroids
$(A_i,[-,-]_i,\zp_{i*})$ associated with principal bundles $P_i(M_i,\zp_i,G_i)$, $i\in\{1,2\}$. It seems
natural to think that $\zF$ induces a vector bundle isomorphism $\zf:A_1\leftrightarrow A_2$ over a
diffeomorphism $\zvf:M_1\leftrightarrow M_2$ that implements $\widetilde{\zF}:{\cal V}_1\leftrightarrow{\cal
V}_2$. However, the preceding guess is not true in general, as shows the example of the Atiyah algebroid
$TM\times\R$ of first-order differential operators on a manifold $M$ \cite{GrabPoncinAuto}. Note that
$TM\times\R$ is isomorphic to the Lie algebroid of linear differential operators acting on smooth sections of
a real vector bundles of rank $1$ over $M$ independently  whether the line bundle is trivial or not
\cite{GrabPoncinDiffOpLineBundles}. We can however build the mentioned vector bundle isomorphism under the
additional assumption that the structure groups $G_i$ are semisimple.\medskip

We are now prepared to state the main result of this section, which yields in particular that the Lie algebra
structure of the space of sections of an Atiyah algebroid recognizes not only the smooth structure of the base
manifold but also the vector bundle structure of the algebroid.

\begin{theo}\label{TheoMain} Let $A_i$, $i\in\{1,2\}$, be smooth Atiyah
algebroids associated with principal bundles with semisimple structure groups. Isomorphisms $\zF:{\cal
A}_1\leftrightarrow{\cal A}_2$ of the Lie algebras ${\cal A}_i$ of smooth sections of the bundles $A_i$ are in
1-to-1 correspondence with isomorphisms $\zf:A_1\leftrightarrow A_2$ of the corresponding Lie algebroids.
\end{theo}

\begin{rem} Eventually, the Lie algebra homomorphism property of $\zF$ might be
encoded in the similar property of the dual vector bundle isomorphism $\zf^*:A^*_2\leftrightarrow A^*_1$ over
$\zvf^{-1}$ (defined in $A^*_{2m}$ by $\zf^*_m:=\,^t\zf_{\zvf^{-1}(m)}$, where notations are self-explaining)
with respect to the linear Poisson structure of $A^*_i$ that is associated with the Lie algebroid structure of
$A_i$.\medskip

Let us recall the basic results pertaining to the mentioned Poisson structure of the dual vector bundle
${\frak p}:L^*\to M$ of any smooth Lie algebroid $(L,[-,-],\widehat{-})$ over a smooth manifold $M$. We set
again ${\cal L}=\zG(L)$. The map $-^{\bullet}:{\cal L}\ni \E\to \E^{\bullet}\in\Ci(L^*)$, which is defined,
for any $\zl\in L^*_m$, $m\in M$, by $\E^{\bullet}(\zl)=\zl(\E_m)\in\R$ and associates to any smooth section
of $L$ a smooth function of $L^*$ that is linear in the fibers, is visibly an injective and nonsurjective
linear mapping. It is well-known, see e.g. \cite{MarlesCalcAldGrdPoisson}, that there is a unique Poisson
structure $\{-,-\}$ on $L^*$, such that $-^{\bullet}$ is a Lie algebra homomorphism, i.e.
$\{\E^{\bullet},\E'^{\bullet}\}=[\E,\E']^{\bullet}$, for any $\E,\E'\in{\cal L}$. Of course, this implies that
$\{-,-\}$ is a linear Poisson bracket. Moreover, it follows from the module-Lie compatibility condition in the
Lie algebroid and the Leibniz property of the Poisson bracket that necessarily $\{\E^{\bullet},g\circ{\frak
p}\}=\widehat{\E}(g)\circ{\frak p}$ and $\{f\circ{\frak p},g\circ{\frak p}\}=0$, for any $f,g\in\Ci(M)$ and
any $\E\in{\cal L}$.
\end{rem}

Roughly speaking, since a fiber $A_{i\,m}$, $m\in M_i$, of $A_i$ can be viewed as the quotient space ${\cal
A}_i/{\cal A}_i(m)$, where ${\cal A}_i(m)=\{a\in{\cal A}_i:a_m=0\}$, any Lie algebra isomorphism $\zF:{\cal
A}_1\leftrightarrow {\cal A}_2$ induces isomorphisms $\zf_m:A_{1m}\leftrightarrow A_{2\zvf(m)}$, where
$\zvf:M_1\leftrightarrow M_2$ is the diffeomorphism generated by Theorem \ref{TheoPursellShanks}, if ${\cal
A}_1(m)$ and ${\cal A}_2(\zvf(m))$ can be characterized in Lie algebraic terms.\medskip

To simplify notations, we drop in the following index $i$. Equation \eqref{AtiyahAldBracketLocal} implies that
the value at $m\in M$ of the bracket $[a,a']$ of two sections $a,a'\in{\cal A}$ is given by \be
[a,a']_m=[a\vert_U,a'\vert_U]_m=([v_1,v_2]_{{\frak X}(U)},[\zg_1,\zg_2]_{\frak
g}+v_1(\zg_2)-v_2(\zg_1))_m,\label{L(m)K(m)Bracket}\ee where $U$ denotes a chart domain of $M$ and a
trivialization domain of the principle bundle $P(M,\zp,G)$, where ${\frak g}=\op{Lie}(G)$, and where
$a\vert_U=(v_1,\zg_1),$ $a'\vert_U=(v_2,\zg_2)\in{\frak X}(U)\times\Ci(U,{\frak g})$. As $a\in{\cal A}(m)$ if
and only if $v_{1m}=\zg_{1m}=0$ and as $a'\in{\cal K}$ entails
$0=(\zp_*a')\vert_U=\zp_*\vert_U(v_2,\zg_2)=v_2$, it follows from Equation \eqref{L(m)K(m)Bracket} that the
${\cal A}(m)$ could be the (maximal) Lie subalgebras of ${\cal A}$, such that $[{\cal A}(m),{\cal K}]\subset
{\cal K}(m)=\{k\in{\cal K}:k_m=0\}$. Hence, the ${\cal A}(m)$ are characterized in Lie algebraic terms, if the
${\cal K}(m)$ are. Equation \eqref{L(m)K(m)Bracket} allows seeing that the ${\cal K}(m)$ are
finite-codimensional Lie ideals of ${\cal K}$. The fact that, for any point $m\in M$ and any (maximal) ideal
${\frak g}_0\subset{\frak g}$, the space ${\cal K}(m,{\frak g}_0):=\{k\in{\cal K}:k_m\in{\frak g}_0\}$ is a
(maximal) finite-codimensional Lie ideal of ${\cal K}$, suggests that ${\cal K}(m)$ is the intersection of all
maximal finite-codimensional Lie ideals ${\cal K}(m,{\frak g}_0)$, as, for a semisimple Lie group $G$, the
intersection of all maximal ideals ${\frak g}_0$ of ${\frak g}$ vanishes.\medskip

These ideas lead to the following theorems. The first one is based upon a series of lemmata.

\begin{rem} Let us stress that in the sequel the structure group
of the principal bundle, which gives rise to the considered Atiyah algebroid, is assumed to be
semisimple.\end{rem}

If $G$ is a semisimple Lie group, its Lie algebra ${\frak g}$ has no nonzero solvable ideals, its radical
${\frak r}$ vanishes and ${\frak g}^{(1)}={\frak g}$. The latter clearly entails ${\cal K}^{(1)}={\cal K}$,
see Equation (\ref{KernelSemisimple}).

\begin{lem}\label{InfDimDerIdeal} The kernel ${\cal K}=\zG(P\times_G{\frak g})$ of the Atiyah
sequence of modules and Lie algebras, which is implemented by a principal bundle $P$ with semisimple structure
group $G$, is an infinite-dimensional Lie algebra, such that ${\cal K}^{(1)}={\cal K}$.\end{lem}

In the following, a ``maximal object with a given property'' is an object that is not strictly contained in
any proper object with the same property.

\begin{lem}\label{ModuleLieIdeals} Let ${\cal K}$ be as detailed in Lemma \ref{InfDimDerIdeal}.
Any maximal finite-codimensional Lie ideal ${\cal K}_0$ in ${\cal K}$ is modular, i.e. ${\cal N K}_0={\cal
K}_0$.\end{lem}

\begin{proof} Since \be[{\cal N K}_0,{\cal K}]=[{\cal K}_0,{\cal N K}]\subset {\cal K}_0\subset {\cal N
K}_0,\label{Modular}\ee the space ${\cal N K}_0$ is a finite-codimensional Lie ideal in ${\cal K}$. Hence,
either ${\cal N K}_0={\cal K}_0$ or ${\cal N K}_0={\cal K}$. In the second case, Equation \eqref{Modular}
shows that ${\cal K}^{(1)}\subset {\cal K}_0$ -- a contradiction in view of the semisimplicity assumption for
$G$ and the maximality assumption for ${\cal K}_0$.\end{proof}

\begin{lem}\label{AssRadicals} Let ${\cal N}=\Ci(M)$ be the associative commutative
unital algebra of smooth functions of a manifold $M$ and let
$I\subset{\cal N}$. Denote by $\op{Spec}({\cal N},I)$ the set of
all maximal finite-codimensional ideals of ${\cal N}$ that contain
$I$ and set $\bar{I}=\cap_{J\in\op{Spec}({\cal N},I)}J$.
Eventually, if $I$ is an ideal in ${\cal N},$ let
$\sqrt{I}=\{f\in{\cal N}:f^n\in I, \mbox{ for some } n>0\}$ be the
radical ideal of $I$. If $I$ is $s$-codimensional in ${\cal N}$,
$s>0$, there are points $m_1,\ldots,m_{\E}\in M$, $1\le\E\le s$,
such that \be \sqrt{I}=\bar
I=\cap_{i=1}^{\E}J(m_i)=\cap_{i=1}^{\E}\{f\in{\cal
N}:f(m_i)=0\}.\ee Moreover, we have $\sqrt{I}^{\,s}\subset I$.
\end{lem}

\begin{proof}

Since the points $m\in M$ are in $1$-to-$1$ correspondence with the maximal finite-codimensional ideals
$J(m)\subset{\cal N}$, see Equation \eqref{Maxfincodimidealsfcts}, we have $\bar I=\cap_{m\in M:J(m)\supset
I}J(m)$. Since each finite-codimensional ideal is included in at least one maximal finite-codimensional ideal,
the preceding intersection is not the intersection of the empty family. Further, as ${\cal
N}/\cap_{i=1}^{\E}J(m_i)\simeq \R^{\E}$ and as $s=\op{dim}{\cal N}/I\ge\op{dim}{\cal N}/\bar I$, it is clear
that $I$ cannot be contained in more than $s$ maximal finite-codimensional ideals, so that $\bar
I=\cap_{i=1}^{\E}J(m_i)$, $\E\in\{1,\ldots, s\}$ and $\op{codim}\bar I=\E\ge 1$.\medskip

If the descending series of finite-codimensional ideals
$$\bar I=I+\bar I\supset I+\bar I^2\supset\ldots\supset I+\bar I^s\supset I+\bar I^{s+1}\supset
I$$ were strictly decreasing (except maybe for the last inclusion), we could not have $\op{codim}\bar I\ge 1$.
Hence, \be I+\bar I^{n}=I+\bar I^{n+1},\label{StabDescIdeals}\ee for some $n\in\{1,\ldots,s\}$. Remember now
Nakayama's lemma that holds true for any commutative ring $R$ with identity $1$, any ideal ${\frak I}$ in $R$,
and any finitely-generated module $M$ over $R$, and which states that if ${\frak I}M=M$, there is $r\in R$,
$r\sim 1$ modulo ${\frak I}$, such that $rM=0$. If ${\frak I}$ is included in the Jacobson radical of $R$,
i.e. in the intersection of all maximal ideals of $R$, then $r$ is invertible and $M=0$. When applying this
upshot to $R={\cal N}/I$, ${\frak I}=\bar I/I$, and $M=(I+\bar I^{n})/I$, where $M$ is actually
finite-dimensional and where ${\frak I}M=M$ in view of Equation \eqref{StabDescIdeals}, we get $\bar
I^n\subset I$, so that $\bar I\subset\sqrt{I}$. The fact that for all $J\in\op{Spec}({\cal N},I)$ the
inclusion $\sqrt I\subset \sqrt J=\sqrt{J(m)}=J(m)=J$ holds true, entails that $\sqrt I\subset \bar I$.
Eventually, $\sqrt I=\bar I$ and $\sqrt I^s\subset \sqrt I^n=\bar I^n\subset I$.
\end{proof}

\begin{theo}\label{LieIdealsCharacterization} Let ${\cal K}$ be the kernel of the Atiyah sequence
of modules and Lie algebras associated to a principal bundle with semisimple structure group $G$ over a
manifold $M$. The maximal finite-codimensional Lie ideals of ${\cal K}$ are exactly the ideals of the form
${\cal K}(m,{\frak g}_0)=\{k\in{\cal K}:k_m\in{\frak g}_0\}$, where $m\in M$ and where ${\frak g}_0$ is a
maximal Lie ideal of ${\frak g}=\op{Lie}(G)$.\end{theo}

\begin{proof} Take a maximal finite-codimensional Lie ideal ${\cal K}_0$ in ${\cal
K}$. The ${\cal N}$-module structure of ${\cal K}_0$, see Lemma \ref{ModuleLieIdeals}, allows switching from
the Lie algebraic to the associative context and then applying Lemma \ref{AssRadicals}. Indeed, in view of
this modularity, the space $I_k:=\{f\in{\cal N}:fk\in{\cal K}_0\}$ is, for any $k\in{\cal K}$, an associative
ideal in ${\cal N}$, which is finite-codimensional as it is the kernel of the linear map $\E:{\cal N}\ni f\to
\widetilde{fk}\in{\cal K}/{\cal K}_0$ with values in a finite-dimensional space. If
$\widetilde{k_1},\ldots,\widetilde{k_q}$ is a basis of the space ${\cal K}/{\cal K}_0$, the intersection
$I=\cap_{i=1}^q I_{k_i}$ is an associative ideal in ${\cal N}$, verifies $I{\cal K}\subset {\cal K}_0$, and is
nonzero- and finite-codimensional ($\op{codim}I\le \sum_{i=1}^q\op{codim}I_{k_i}$). Hence, Lemma
\ref{AssRadicals} is valid for this ideal $I$; in the sequel, we use the notations of this lemma.\medskip

It is obvious that the space $\sqrt{I}{\cal K}$, which is made up
by finite sums of products $fk$, $f\in\sqrt{I},$ $k\in{\cal K}$,
is a Lie ideal in ${\cal K}$ and consists of the sections in $\cK$
which vanish at the points $m_1,\ldots,m_{\E}$ that provide the
radical $\sqrt I$. Indeed, if we take a partition of unity
$(U_{\frak a},\zg_{\frak a})_{\frak a\in\frak A}$ that is
subordinated to a cover by domains of local coordinates, such a
section is locally, in each $U_{\frak a}$, a combination of local
sections and local functions that vanish at the $m_i\in U_{\frak
a}$. It suffices then to use the partition $\zg_{\frak a}$ to show
that the considered section is actually an element of $\sqrt I
\cK$.\medskip

This characterization of the elements of the Lie ideal $\sqrt I{\cal K}$ in ${\cal K}$ allows proving that
there exists a class in the Lie algebra ${\cal K}/\sqrt I{\cal K}$ that does not contain any element of ${\cal
K}_0$.

Indeed, if any class contained an element of ${\cal K}_0$, then,
for any $k\in{\cal K}$, we would have a series of equations
$$\begin{array}{c}k=k_{0}+\sum_{i_1}r_{i_1}k_{i_1},\\ k_{i_1}=k_{0i_1}+\sum_{i_2}r_{i_1i_2}k_{i_1i_2},\\
\ldots\\ k_{i_1\ldots i_{s-1}}=k_{0i_1\ldots
i_{s-1}}+\sum_{i_{s}}r_{i_1\ldots i_{s}}k_{i_1\ldots
i_{s}},\end{array}$$ with $k_{0i_1\ldots i_{u-1}}\in{\cal K}_0$,
$r_{i_1\ldots i_u}\in\sqrt I$, and $k_{i_1\ldots i_u}\in{\cal K}$.
Thus,
$$k=\lp k_{0}+\sum_{i_1}r_{i_1}k_{0i_1}+\ldots +\sum_{i_1\ldots i_{s-1}}r_{i_1}\ldots r_{i_1\dots i_{s-1}}k_{0i_1\ldots i_{s-1}}\rp
+ \sum_{i_1\ldots i_{s}}r_{i_1}\ldots r_{i_1\ldots
i_{s}}k_{i_1\ldots i_{s}},$$ and, since ${\cal K}_0$ is modular in
view of Lemma \ref{ModuleLieIdeals}, the parenthesis in the
{\small RHS} of this equation is an element of ${\cal K}_0$,
whereas the last term is in ${\cal K}_0$ as well, since $\sqrt
I^s\subset I$, due to Lemma \ref{AssRadicals}, and $I{\cal
K}\subset {\cal K}_0$. It follows that ${\cal K}\subset {\cal
K}_0$ - a contradiction, since ${\cal K}_0$ is maximal by
assumption.\medskip

We just proved that there exists at least one element $k\in{\cal K}$, such that for any $k_0\in{\cal K}_0$, we
have $k-k_0\notin\sqrt I{\cal K}$. But then, there is at least one point $m\in M$, such that the space ${\frak
g}_0:=\{k_{0m}:k_0\in{\cal K}_0\}$ is a proper subspace of ${\frak g}$.

In fact, otherwise ${\frak g}_0$ would vanish for any $m\in M$ or would coincide with ${\frak g}$ for any
$m\in M$. As the first alternative is impossible, it follows that for any $m\in M$ and any ${\zk}\in{\cal K}$,
there exists ${\zk}_0\in{\cal K}_0$, such that $\zk_m=\zk_{0m}$. This assertion holds true in particular for
the points $m_1,\ldots,m_{\E}\in M$ and the above section $k\in{\cal K}$; there are $k_{0i}\in{\cal K}_0$,
such that $k_{m_i}=k_{0i,m_i}$. Take now open subsets $U_i\subset M$ that contain $m_i$ and have pairwise
empty intersections, as well as bump functions $\za_i\in{\cal N}$ with value $1$ in a neighborhood of $m_i$
and compact support in $U_i$, and set finally $k_0=\sum_i\za_ik_{0i}\in{\cal K}_0.$ It is clear that
$k_{0,m_j}=k_{m_j}$, so that $k-k_0\in{\cal K}$ vanishes at $m_1,\ldots,m_{\E}$ and is thus an element of
$\sqrt I{\cal K}$. Since this is impossible, there actually exists at least one point $m\in M$, such that
${\frak g}_0:=\{k_{0m}:k_0\in{\cal K}_0\}$ is proper in ${\frak g}$.\medskip

It is readily checked that the subspace ${\frak g}_0$ is a Lie ideal in ${\frak g}$. Note further that for any
point $p\in M$ and any Lie ideal ${\frak h}_0\subset {\frak g}$, the subspace ${\cal K}(p,\frak
h_0):=\{k\in{\cal K}:k_p\in\frak h_0\}$ is a Lie ideal in ${\cal K}$, and that this ideal is
finite-codimensional since the space ${\cal K}/{\cal K}(p,\frak h_0)$ is isomorphic to ${\frak g}/{\frak
h_0}$. As obviously ${\cal K}_0\subset{\cal K}(m,{\frak g}_0),$ the preceding observation entails that the
maximal finite-codimensional ideal ${\cal K}_0$ is included in the finite-codimensional ideal ${\cal
K}(m,\frak g_0)$, which is in turn strictly included in ${\cal K}$, since $\frak g_0$ is proper in $\frak g$.
Hence, ${\cal K}_0={\cal K}(m,\frak g_0)$. Eventually, the ideal $\frak g_0$ is maximal in $\frak g$;
otherwise, it would be strictly included in a proper ideal $\frak h_0$, so that ${\cal K}_0={\cal K}(m,\frak
g_0)$ would on his part be strictly included in the proper finite-codimensional ideal ${\cal K}(m,\frak h_0)$.
This concludes the proof of the first part of Theorem \ref{LieIdealsCharacterization}.\medskip

Let now $m\in M$ and let $\frak g_0$ be a maximal Lie ideal in $\frak g$. We already mentioned that ${\cal
K}(m,\frak g_0)$ is then a finite-codimensional Lie ideal in ${\cal K}$. If this ideal is not maximal it is
strictly included in a proper finite-codimensional ideal ${\cal K}_1$, which can of course be assumed to be
maximal. But in this case, the first part of the theorem implies that ${\cal K}_1={\cal K}(p,\frak h_0)$, for
some $p\in M$ and some maximal ideal $\frak h_0$ in $\frak g$. Hence, $k_m\in\frak g_0$ entails $k_p\in\frak
h_0$, for any $k\in{\cal K}$. From this it first follows that $m=p$. Indeed, otherwise we take $k\in{\cal K}$,
such that $k_p\notin\frak h_0$, as well as a bump function $\za\in\Ci(M)$ that has value $0$ in a neighborhood
of $m$ and value $1$ in some neighborhood of $p$. The fact that $(\za k)_m\in\frak g_0$ and $(\za
k)_p\notin\frak h_0$ then constitutes a contradiction. We now see that the maximal ideal $\frak g_0$ is
included in the (maximal and thus) proper ideal $\frak h_0.$ This shows that $\frak g_0=\frak h_0$ and, since
$m=p,$ that ${\cal K}(m,\frak g_0)={\cal K}(p,\frak h_0)={\cal K}_1$, so that ${\cal K}(m,\frak g_0)$ is
actually maximal.
\end{proof}

The space ${\cal K}(m,\frak g_0)$ is a maximal finite-codimensional Lie ideal in the Lie ideal ${\cal K}$ in
${\cal A}$. The next proposition provides the normalizer $N({\cal A},{\cal K}(m,\frak g_0))$ of ${\cal
K}(m,\frak g_0)$ in ${\cal A}$, i.e. the biggest Lie subalgebra of ${\cal A}$ that admits ${\cal K}(m,\frak
g_0)$ as a Lie ideal.

\begin{prop} The normalizer in ${\cal A}$ of any maximal finite-codimensional
Lie ideal ${\cal K}(m,\frak g_0)$ in ${\cal K}$ coincides with the maximal finite-codimensional Lie subalgebra
${\cal A}_m=\{a\in{\cal A}:\widehat a_m=0\}$.\end{prop}

\begin{proof} Let us look for all the $a\in{\cal A}$, such that $[a,k_0]\in{\cal K}(m,\frak
g_0)$, for any $k_0\in{\cal K}(m,\frak g_0)$. Since this condition only involves the value of the local
bracket $[-,-]$ at the point $m$, we consider a chart domain $U$ of $M$ around $m$ that is simultaneously a
trivialization domain of principal bundle that gives rise to the Atiyah sequence; we choose the local
coordinates $x=(x_1,\ldots,x_n)$ in $U$ in such a way that $x(m)=0$. If we set $a\vert_U=(v,\zg)$ and
$k_0\vert_U=(0,\zh)$, the condition reads
$$[\zg_m,\zh_m]_{\frak g}+v_m(\zh)\vert_m\in\frak g_0,$$ for any $\zh\in\Ci(U,\frak
g)$ such that $\zh_m\in\frak g_0$. Since the bracket in the first term is always in $\frak g_0$, the condition
means that $v_m=0$, as we easily see when taking $\zh=x_ke$, where $e\in\frak g\setminus\frak g_0$ is a
nonzero vector. Eventually, the normalizer is exactly ${\cal A}_m=\{a\in{\cal A}:\widehat a_m=0\}$.\end{proof}

\begin{theo}\label{LieIdealsIntersection} Let $0\to {\cal K}\to {\cal A}\to {\cal V}\to 0$
be the Atiyah sequence of modules an Lie algebras associated to a principal bundle with semisimple structure
group over a manifold $M$. The intersections of all maximal finite-codimensional Lie ideals in ${\cal K}$,
which are normalized by a given maximal finite-codimensional Lie subalgebra in ${\cal A}$, are exactly the
${\cal K}(m)=\{k\in{\cal K}:k_m=0\}$, $m\in M$, where $m$ is the point that characterizes the chosen
subalgebra of ${\cal A}$.\end{theo}

\begin{proof} There is a unique $m\in M$ such that ${\cal A}_m$ coincides with the
chosen maximal finite-codimensional Lie subalgebra. In view of the preceding proposition, the corresponding
normalized maximal finite-codimensional Lie ideals of ${\cal K}$ are exactly the ${\cal K}(m,\frak g_0)$,
where $\frak g_0$ runs through the maximal ideals of $\frak g$. Hence, the wanted intersection is ${\cal
K}(m):=\{k\in{\cal K}:k_m=0\}$, as the intersection of all maximal ideals vanishes in a semisimple Lie
algebra. \end{proof}

\begin{theo}\label{AmLieCharact} Consider an Atiyah sequence as in Theorem
\ref{LieIdealsIntersection}. The Lie subalgebras ${\cal A}(m):=\{a\in{\cal A}:a_m=0\}$, $m\in M$, of ${\cal
A}$ can be characterized as the maximal Lie subalgebras ${\cal A}_0$ in ${\cal A}$ such that $[{\cal
A}_0,{\cal K}]$ is included in a certain ${\cal K}(m)$, $m\in M$. In other words, the Lie subalgebras ${\cal
A}(m)$, $m\in M$, are exactly the maximal Lie subalgebras ${\cal A}_0$ in ${\cal A}$ such that $[{\cal
A}_0,{\cal K}]$ is included in the intersection of all maximal finite-codimensional Lie ideals in ${\cal K}$,
which are normalized by a given maximal finite-codimensional Lie subalgebra in ${\cal A}$. The maximal
finite-codimensional Lie subalgebra that corresponds to a precise ${\cal A}(m)$ is ${\cal A}_m$.
\end{theo}

\begin{rem} Let us emphasize that in this statement ``maximal''
means ``maximal in the class of all Lie subalgebras ${\cal A}_0$ of ${\cal A}$ that have the property $[{\cal
A}_0,{\cal K}]\subset{\cal K}(m)$, for some $m\in M$''.\end{rem}

\begin{proof} Let $a\in{\cal A}$, $k\in{\cal K}$, and $m\in M$. Consider again a chart and trivialization domain
$U$ around $m$ and set $a\vert_U=(v,\zg)$ and $k\vert_U=(0,\zh)$. Then, \be [a,k]_m=[\zg_m,\zh_m]_{\frak
g}+v_m(\zh)\vert_m\label{FinalLieCharTheoBracket}.\ee

Equation \eqref{FinalLieCharTheoBracket} entails that any Lie subalgebra ${\cal A}(m)$ verifies $[{\cal
A}(m),{\cal K}]\subset {\cal K}(m)$.\medskip

Conversely, let ${\cal A}_0$ be any maximal Lie subalgebra of ${\cal A}$ such that $[{\cal A}_0,{\cal
K}]\subset {\cal K}(m)$, for a certain $m\in M$. It follows from Equation \eqref{FinalLieCharTheoBracket},
written for $\zh\in\frak g\subset\Ci(U,\frak g)$, that any $a\in {\cal A}_0$ has a vanishing second component
$\zg_m$, since the center of a semisimple Lie algebra vanishes. When writing now Equation
\eqref{FinalLieCharTheoBracket} for an arbitrary function $\zh\in\Ci(U,\frak g)$, we get $v_m(\zh)\vert_m=0$.
Thus $v_m=0$ and ${\cal A}_0\subset{\cal A}(m)$. If ${\cal A}_0\subsetneqq{\cal A}(m)$, we have ${\cal
A}_0\subsetneqq{\cal A}(m)\subset{\cal A}_m\subsetneqq{\cal A}$, since ${\cal A}_m$ is maximal. As ${\cal
A}_0$ is maximal, it follows that ${\cal A}_0={\cal A}(m).$\medskip

It is now easily seen that ${\cal A}(m)$ is maximal among the Lie subalgebras ${\cal A}_0$ of ${\cal A}$ such
that $[{\cal A}_0,{\cal K}]\subset {\cal K}(p)$, for some $p\in M$. Indeed, since ${\cal A}/{\cal A}(m)\simeq
A_m$, where $A_m$ denotes the fiber at $m$ of the Atiyah algebroid $A$, the space ${\cal A}(m)$ is of course
finite-codimensional. So if ${\cal A}(m)$ is strictly included in a proper Lie subalgebra ${\cal A}_0$ such
that $[{\cal A}_0,{\cal K}]\subset{\cal K}(p)$, $p\in M,$ we can assume that ${\cal A}_0$ is maximal in the
considered class. But then ${\cal A}_0={\cal A}(p)$, $p\in M$, and ${\cal A}(m)\subsetneqq{\cal A}(p)$. The
usual argument based upon a smooth function that has value $0$ at $m$ and value $1$ at $p$ then shows that
$m=p$ and that ${\cal A}(m)$ is maximal.\end{proof}

We are now able to provide the proof of Theorem \ref{TheoMain}.

\begin{proof}[Proof of Theorem \ref{TheoMain}] Let $\zF:{\cal A}_1\leftrightarrow {\cal
A}_2$ be a Lie algebra isomorphism.\medskip

Since Theorem \ref{AmLieCharact} characterizes the subalgebras ${\cal A}(m)$, $m\in M$, in pure Lie algebraic
terms, isomorphism $\zF$ transforms an ${\cal A}_1(m)$ into an ${\cal A}_2(p)$. Therefore, ${\cal
A}_2(p)=\zF({\cal A}_1(m))\subset\zF({\cal A}_{1\,m})={\cal A}_{2\,\zvf(m)},$ where $\zvf$ is the
diffeomorphism that is implemented by $\zF$, see Theorem \ref{TheoPursellShanks}. But then $p=\zvf(m)$.
Indeed, let $v\in{\cal V}_2$ be a vector field of $M_2$ such that $v_{\zvf(m)}\neq 0$ and take $a\in{\cal
A}_2$ with anchor $\widehat a=v.$ If $p\neq \zvf(m)$, we can choose a function $\za$ with value $0$ around $p$
and value $1$ around $\zvf(m)$. Thus $\za a\in{\cal A}_2(p)\subset{\cal A}_{2\,\zvf(m)}$, so that $\widehat
a_{\zvf(m)}=v_{\zvf(m)}=0$. It follows that
\be\label{TheoMainObservationMain}\zF({\cal A}_1(m))={\cal
A}_2(\zvf(m)).\ee

As the fiber $A_{i\,m}$ is isomorphic to the vector space ${\cal
A}_i/{\cal A}_i(m)$, the preceding upshot entails that $\zF$
induces linear maps $\zf_m:A_{1\,m}\to A_{2\,\zvf(m)}$, as well as
a smooth vector bundle map $\zf:A_1\to A_2$ over the
diffeomorphism $\zvf:M_1\leftrightarrow M_2.$ Smoothness of $\zf$
is a consequence of the fact that the map $\zF$, which $\zf$
induces between sections, transforms smooth sections into smooth
sections. The bundle map $\zf$ is actually a vector bundle
isomorphism over $\zvf$, since $\zf_m$ is bijective, due to
bijectivity of $\zF$ and Equation
(\ref{TheoMainObservationMain}).\end{proof}

To conclude we combine the preceding upshots with results of \cite{JKub} and get the following
characterization of isomorphisms of Lie algebras of semisimple Atiyah algebroids.

\begin{theo}\label{TheoMain2} Let $A_i$, $i\in\{1,2\}$, be smooth
Atiyah algebroids associated with principal bundles $P_i(M_i,\zp_i,G_i)$ with semisimple structure groups. The
Lie algebras $\cA_i$ of smooth sections of the bundles $A_i$ are isomorphic if and only if the Lie algebroids
$A_i$ are isomorphic, or, as well, if and only if the principal bundles $P_i$ are locally
isomorphic.\end{theo}



\section{Isomorphisms of Lie algebras of Atiyah algebroids -
reductive structure groups}

Let \be 0\to K\stackrel{i}{\to}A\stackrel{\zp_*}{\to}TM\to 0\label{AthiAlg}\ee be an Atiyah algebroid
associated to a principal bundle $P(M,\zp,G)$ with connected {\it reductive structure group} $G$. The Lie
algebra $\frak g=\op{Lie}(G)$ then canonically splits into its Abelian center $Z\frak g$ and its semisimple
derived ideal $[\frak g,\frak g]$. It follows from Equation (\ref{KernelReducive}) that the kernel $K$ is
similarly split,
$$K=\sZ\oplus [K,K],$$
where we wrote simply $\sZ$ for $ZK$, and from Equations (\ref{SectionsCenter}) and (\ref{SectionsDerived}) that
$${\cal K}=\zG(K)=\zG(\sZ)\oplus\zG([K,K])={\cZ}\oplus [{\cal K},{\cal K}].$$

\begin{prop}\label{TrivialCenter} The center sub-{\small LAB} $\sZ=P\times_G Z\frak g$ of the kernel {\small LAB} $K=P\times_G\frak g$ of
the Atiyah sequence of a principal bundle $P(M,\zp,G)$ with connected structure Lie group $G$ is the trivial
bundle $\sZ=M\times Z\frak g$ and it admits a global frame made up by constant functions $\frak
c_1,\ldots,\frak c_k\in\Ci(M,Z\frak g)\simeq\Ci(M,\R^k)$.\end{prop}

\begin{proof} Let us first recall that the adjoint action of a
connected finite-dimensional Lie group $G$ on the center $Z\frak g$ of its Lie algebra is trivial.

Triviality of the adjoint action on the kernel entails that the sections of $\sZ=P\times_GZ\frak g$, i.e. the
$G$-invariant functions from $P$ to $Z\frak g$, are exactly the functions from $M$ to $Z\frak g$. Hence, any
basis $\frak z_1,\ldots,\frak z_k$ of the vector space $Z\frak g$ corresponds to global sections $\frak
c_1,\ldots,\frak c_k$ of $\sZ$ that obviously form a global frame.
\end{proof}

Since $\sZ$ is a Lie algebroid ideal of $A$, the sequence \be 0\to
\sZ\stackrel{i}{\to}A\stackrel{p}{\to}A/\sZ\to 0,\label{QuotiSeq}\ee is a short exact sequence of Lie
algebroids and the transitive quotient Lie algebroid $\wA:=A/\sZ$ is associated with the short exact sequence
\be 0\to K/\sZ\simeq [K,K]\stackrel{\tilde \imath}{\to}
\wA\stackrel{\tilde \zp_*}{\to} TM\to 0,\label{SesQuotiLad}\ee see Section \ref{Lads}. Moreover, the isotropy
algebra of $\wA$, the quotient algebra $\widetilde{\frak g}:={\frak g}/Z\frak g$, is semisimple. We will write
$\wcA:=\zG(A/\sZ)={\cal A}/{\cZ}$ for the Lie algebra of sections of this Lie algebroid and $\hat{a}$ for the
anchor $\tilde \zp_*(\tilde a)$ of a section $\tilde a\in\wcA$.

\begin{rem} It is known \cite{JKub} that any transitive Lie algebroid with semisimple
{\small LAB} is the Atiyah algebroid of some principal bundle. Hence, the quotient Lie algebroid $\wA=A/\sZ$ is an
Atiyah algebroid, namely that of the principal bundle $P/ZG(M,\tilde \zp,G/ZG)$, where notations are
self-explaining and where $G/ZG$ is semisimple.
\end{rem}

In the sequel, we use the global frame $\frak c_1,\ldots,\frak c_k\in\Ci(M,Z\frak g)$ of ${\cZ}$ made up by
constant functions.

\begin{lem}\label{CoincActions} For any sections $a\in{\cal A}$ and $c=\sum_if^i\frak c_i\in {\cZ}$, $f^i\in\Ci(M)$,
the adjoint action on $c$ by $a$ and the canonical action by $\zp_*a\in{\frak X}(M)$ coincide:
$$[a,c]=(\zp_*a)(c)=\sum_i(\zp_*a)(f^i)\frak c_i.$$\end{lem}

\begin{proof} Since $\frak c_i$ is a constant $Z\frak g$-valued function, (it follows for instance
from the local form of the Atiyah algebroid bracket that) we have $[a,\frak c_i]=0$. The lemma is then an
immediate consequence of the Leibniz property of the Lie bracket $[-,-]$.\end{proof}

\begin{prop}\label{LaCharZKm} The set ${\cZ}(m)$, $m\in M$, of all the sections in
${\cZ}$ that vanish at $m$ is given in Lie algebraic terms by
$${\cZ}(m)=[{\cal A}_m,{\cZ}],$$ where ${\cal A}_m=\{a\in{\cal
A}:(\zp_*a)_m=0\}$.\end{prop}

\begin{proof} Lemma \ref{CoincActions} entails that the inclusion $\supset$ holds true.
As for the converse inclusion, remember first that any function of $M$ can be written as a sum of Lie
derivatives \cite{GrabowskiIsomIdeal}; in particular, there are vector fields $X_j\in{\frak X}(M)$ and
functions $g^j\in\Ci(M)$ such that $1=\sum_jL_{X_j}g^j.$ Let now $c=\sum_i f^i \frak c_i\in {\cZ}(m)$, so that
$f^i(m)=0$, and set $Y_j^i=f^iX_j\in{\frak X}(M)$, so that $f^i=\sum_jL_{Y^i_j}g^j$ and $Y_{j,m}^i=0$. If
$a^i_j\in{\cal A}$ denotes any preimage of $Y_j^i$ by $\zp_*$, we have
$$c=\sum_{ij}(\zp_*a^i_j)(g^j)\frak c_i=\sum_{ij}[a^i_j,g^j\frak c_i]\in [{\cal A}_m,{\cZ}].$$\end{proof}

In the following, we investigate isomorphisms $$\zF:{\cal A}_1\leftrightarrow {\cal A}_2$$ of Lie algebras
${\cal A}_i$ of Atiyah algebroids $A_i$ associated to principal bundles $P_i(M_i,\pi_i,G_i)$ with connected
{\it reductive structure groups} $G_i$. Let $0\to K_i\to A_i\to TM_i\to 0$ be the corresponding Atiyah
sequences. In view of Theorem \ref{TheoPursellShanks}, we have $\Phi(\cK_1)=\cK_2$, so that moreover
$\Phi(\cZ_1)=\cZ_2$. Hence, the isomorphism $\Phi:{\cal A}_1\leftrightarrow{\cal A}_2$ induces an isomorphism
$$\Phi^0:\cZ_1\leftrightarrow \cZ_2$$ of the centers $\cZ_i$
and an isomorphism
$${\Phi}^s:\wcA_1\leftrightarrow\wcA_2$$ of the Lie
algebras $\wcA_i=\zG(\wA_i)$ of the Atiyah algebroids $\wA_i=A_i/\sZ_i$ implemented by principal bundles
with the semisimple structure groups $G_i/ZG_i$. Note that, in view of Theorem \ref{TheoMain2}, the Lie
algebra isomorphism $\zF^s$ is implemented by a Lie algebroid isomorphism $\zf^s:\wA_1\leftrightarrow
\wA_2$ covering a diffeomorphism ${\zvf}:M_1\leftrightarrow M_2$.

\begin{prop}\label{Phi0Prop} The Abelian Lie algebra isomorphism $\zF^0:\cZ_1\leftrightarrow \cZ_2$
is implemented by a vector bundle isomorphism
$\zf^0:\sZ_1\leftrightarrow \sZ_2$ that covers ${\zvf}$. Moreover,
$\zF^0$ is, for any $c_1=\sum_i f^i\frak c_{1i}\in \cZ_1,$ given
by \be \zF^0\lp\sum_i f^i\frak c_{1i}\rp=\sum_{ij} I^i_j \lp
f^j\circ \zvf^{-1}\rp\frak c_{2i},\label{Phi0Equation}\ee where
$I=(I^i_j)\in\op{GL}(k,\R)$.\end{prop}

\begin{proof} Since, due to Proposition \ref{LaCharZKm}, the isomorphism $\zF^0$ generates a bijection
between the sets $\{\cZ_1(m_1):m_1\in M_1\}$ and $\{\cZ_2(m_2):m_2\in M_2\}$, it is implemented by a vector
bundle isomorphism $\zf^0:\sZ_1\leftrightarrow \sZ_2$ between the trivial center bundles, which covers a
diffeomorphism $\zvf^0:M_1\leftrightarrow M_2$:
$$\zf^0:\sZ_1\ni(m_1,z_1)\leftrightarrow (\zvf^0(m_1),{\cal I}(m_1)z_1)\in
\sZ_2,$$ where $\cI(m_1)$ is a vector space isomorphism from $Z\frak g_1$ onto $Z\frak g_2.$ Therefore,
\be\zF^0\lp\sum_if^i\frak c_{1i}\rp=\zf^0\lp \sum_i
f^i(\zvf^{0})^{-1}\frak c_{1i}\rp=\sum_{ij}\; I_i^j(\zvf^0)^{-1}\cdot f^i(\zvf^{0})^{-1}\;\frak
c_{2j},\label{Phi0_1}\ee where compositions are understood and where $I\in\Ci(M_1,\op{GL}(k,\R))$ denotes the
matrix of $\cI$ in the bases $(\frak c_{1i})$ and $({\frak c_{2i}})$.

Let $a_1\in\cA_1$, $c_1\in \cZ_1$ and set $\zp_{1*}a_1=:X_1$. The Lie algebra morphism property of $\zF$ then
reads $$ \zF^0(L_{X_1}c_1)=\zF((\zp_{1*}a_1)(c_1))=\zF[a_1,c_1]=[\zF a_1,\zF^0c_1]=(\zp_{2*}\zF
a_1)(\zF^0c_1).$$ As Theorem \ref{TheoPursellShanks} implies that $\zF^s$ induces a Lie algebra isomorphism
$\widetilde\zF^s={\zvf}_*$ between the algebras of vector fields, we get
$$\zp_{2*}\zF a_1=\tilde\zp_{2*}p_2\zF
a_1=\tilde\zp_{2*}\zF^sp_1a_1=\widetilde\zF^s(\tilde\zp_{1*}p_1a_1)={\zvf}_*\zp_{1*}a_1={\zvf}_*X_1.$$ The
combination of the last two upshots finally gives \be \zF^0(L_{X_1}c_1)=
L_{{\zvf}_*X_1}(\zF^0c_1).\label{IsomProp_ac}
\ee

To simplify notations, denote by $f$ the $\R^k$-valued function with components $f^i$. The combination of the
equations (\ref{IsomProp_ac}) and (\ref{Phi0_1}) then leads to
\be\label{ConstIPhisPhi0}L_{{\zvf}_*X_1}\lp I(\zvf^0)^{-1}\cdot
f(\zvf^{0})^{-1}\rp= I(\zvf^0)^{-1}\cdot (L_{X_1} f)(\zvf^{0})^{-1}= I(\zvf^0)^{-1}\cdot \zvf^{0}_* L_{X_1}
f.\ee For any vector field $X_1\in\frak X(M_1)$ that vanishes at an arbitrarily chosen point $m_1\in M_1$, if
we write Equation (\ref{ConstIPhisPhi0}) at ${\zvf}m_1$, we get $(L_{X_1}f)((\zvf^0)^{-1}{\zvf}m_1)=0$, for
any $f\in\Ci(M_1,\R^k)$. It follows that $X_1$ vanishes at $(\zvf^0)^{-1}{\zvf}m_1,$ if, as assumed, $X_1$
vanishes at $m_1$. Hence, ${\zvf}=\zvf^0$.

Equation (\ref{ConstIPhisPhi0}) gives now
$$L_{{\zvf}_*X_1}\lp I{\zvf}^{-1}\cdot {\zvf}_*f\rp=I{\zvf}^{-1}\cdot L_{{\zvf}_*X_1}
{\zvf}_*f,$$ so that the matrix $I$ is actually constant.
\end{proof}

We now aim at writing $\zF:\cA_1\leftrightarrow\cA_2$ by means of $\zF^0:\cZ_1\leftrightarrow \cZ_2$ and
$\zF^s:\wcA_1\leftrightarrow\wcA_2$. This requires the use of a connection.

\begin{lem} Any right splitting $\n$ of the vector bundle sequence (\ref{AthiAlg}),
i.e. any connection of an Atiyah algebroid $A$ with connected reductive structure group, naturally induces a
right splitting $\tilde \n$ of the sequence (\ref{SesQuotiLad}) and a right splitting $\bar\n$ of the sequence
(\ref{QuotiSeq}). Moreover, for any $k^{(1)}\in[K,K]$, we have $\bar\n p(k^{(1)})=k^{(1)}$. The preceding
splitting allows identifying $\wA$ with a vector subbundle $\bar\n (\wA)$ of $A$ that verifies
$A=\sZ\oplus \bar\n (\wA)$ (as vector bundles) and $[K,K]\subset \bar\n (\wA)$.\end{lem}

\begin{proof} The first claim is clear, it suffices to set $\tilde\n=p\circ\n$.
As for $\bar\n$, observe that, if $'$ (resp. $''$) is the projection of $K=\sZ\oplus [K,K]$ onto $\sZ$ (resp.
$[K,K]$) and if $p(a)\in \wA$, the difference $a-\n\zp_*a$ belongs to $K$ and the sum
$\n\zp_*a+(a-\n\zp_*a)''$ is well-defined in $\wA$. The bundle map
$$\bar\n:\wA\ni p(a)\to \n\zp_*a+(a-\n\zp_*a)''\in A$$ is then the
searched splitting $\bar\n$, since
$$p(\bar\n
p(a))=p(\n\zp_*a+(a-\n\zp_*a)'')=p(\n\zp_*a+(a-\n\zp_*a)'+(a-\n\zp_*a)'')=p(a).$$ The assertion $\bar\n\circ
p=\op{id}$ on {[K,K]} immediately follows from the definition of $\bar\n$, whereas the last part of the lemma
is obvious.\end{proof}

\begin{lem} Let $\n$ be any connection of an Atiyah
algebroid $A\to M$ with connected reductive structure group. If we transfer the Lie algebroid bracket on
${\cal A}={\cZ}\oplus\bar\n({\wcA})$ from $\cA$ to ${\cZ}\oplus{\wcA}$, it reads, for any $c,c'\in {\cZ}$ and
$\wa,\wa'\in\wcA$, \be [\![c+\wa,c'+\wa']\!]_{\zw}=\hat a(c')-\hat a'(c)+\omega(\hat a,\hat
a')+[\wa,\wa']^{\;\tilde{}}, \label{LadBracketThroughIsom}\ee where, for any $X,Y\in{\frak X}(M)$,
\be\label{omega}\omega(X,Y)=\sum_j\omega^j(X,Y)\frak c_j,\ee for some closed 2-forms $\omega^j$ on $M$. In
other words, any Atiyah algebroid $(A,[-,-],\zp_*)$ over $M$ with reductive structure group is isomorphic with
a {\rm model reductive Atiyah algebroid} $(\sZ\oplus_M\wA,[\![-,-]\!]_\zw,\wzp_*^0)$, where $\sZ$ is a trivial
bundle and $\wA$ is an Atiyah algebroid over $M$ with semisimple structure group, with the Lie bracket of
sections (\ref{LadBracketThroughIsom}) associated with a closed 2-form $\zw$ on $M$ with values in $Z\frak g$,
and with the anchor map $\wzp_*^0(c+\tilde a)=\wzp_*(\tilde a)=\hat a$.
\end{lem}

\begin{proof} It is well-known that the curvature $R_{\bar\n}$ of $\bar\n$ is a
closed Lie algebroid $2$-form of $\wA$ valued in $\sZ$, so that
$$\zW_{\n}:=R_{\bar\n}\in\zG(\w^2\wA^*\otimes \sZ)=\zG(\w^2\wA^*\otimes Z\frak
g),$$ see e.g. Section \ref{Lads}. To get the transferred bracket, note that the first term $[c,c']$, $c,c'\in
{\cZ}$, of $[c+\bar\n p(a),c'+\bar\n p(a')]$ vanishes, that the second and third terms are of the type
$$[\bar\n p(a),c']=(\zp_*\bar\n p(a))(c')=(\tilde\zp_*p(a))(c')=(\zp_*a)(c')\in {\cZ},$$ due to Lemma
\ref{CoincActions}, whereas the fourth term reads
$$[\bar\n p(a),\bar\n
p(a')]=\zW_{\n}(p(a),p(a'))+\bar\n[p(a),p(a')]^{\;\tilde{}}\in {\cZ}\oplus\bar\n ({\wcA}).$$ Hence the
announced result up to the third term of the {\small RHS} of Equation (\ref{LadBracketThroughIsom}).

We can conclude that the curvature $\zW_{\n}$ is defined on $\frak
X(M)\simeq\wcA / \wcK$, where $\wcK:=\cK/\cZ$, if we prove that it
vanishes once one of the arguments is in
$[\cK,\cK]=\cK^{(1)}\simeq p(\cK^{(1)})=\wcK$. However, since
$\bar\n\circ p=\op{id}$ on $[K,K]$ and since $[K,K]$ is a Lie
algebroid ideal in $A$, we have \be\label{VanishingzW}
\zW_{\n}(p(a),p(k^{(1)}))=[\bar\n p(a),k^{(1)}]-\bar\n
p[a,k^{(1)}]=[\bar\n p(a)-a,k^{(1)}]=0,\ee where the last member
vanishes since by definition $\bar\n p(a)-a=(\n\zp_*a-a)'\in \cZ$.
The resulting $\sZ$-valued (i.e. $Z\frak g$-valued or
$\R^k$-valued) 2-form $\zw=\zw_{\n}$ of $M$ is still closed.
Indeed, this form can be computed, for any $X,Y\in\frak X(M)$, by
$\zw_{\n}(X,Y)=\zW_{\n}(\tilde\n X,\tilde\n Y)$, since
$\tilde\zp_*\tilde\n X=X$. The de Rham differential
$(d\zw_{\n})(X,Y,Z)$ is thus made up by two types of terms,
$$X.\zw_{\n}(Y,Z)=(\tilde\zp_*\tilde\n X).\zW_{\n}(\tilde\n Y,\tilde\n Z)$$
and
$$\zw_{\n}([X,Y],Z)=\zW_{\n}(\tilde\n[X,Y],\tilde\n Z)+\zW_{\n}(R_{\tilde\n}(X,Y),\tilde\n Z)=\zW_{\n}
([\tilde\n X,\tilde\n Y]^{\;\tilde{}},\tilde\n Z),$$
where the first equality is due to Equation (\ref{VanishingzW}) and to the fact $R_{\tilde\n}(X,Y)\in
p(\cK^{(1)})$. Eventually, the considered de Rham differential of $\zw_{\n}$ coincides with the Lie algebroid
differential $(d\zW_{\n})(\tilde\n X,\tilde\n Y,$ $\tilde\n Z)=0$ of the closed form $\zW_{\n}$ of the Lie
algebroid $(\wA,[-,-]^{\;\tilde{}},\tilde\zp_*)$.\end{proof}

Just as the Lie bracket, the Lie isomorphism
$\zF:\cA_1\leftrightarrow\cA_2$ catches a twist when read through
the isomorphism \be\cA_i=\cZ_i\oplus\bar\n_i(\wcA_i)\simeq
\cZ_i\oplus\wcA_i.\label{IsomConnection}\ee Before formulating our
main result, let us recall that on every manifold $M$ there exists
a {\it divergence}, i.e. a linear operator $\op{div}:{\frak
X}(M)\rightarrow C^\infty(M)$, which is a cocycle, i.e.
$$\op{div}[X,Y]=X\left(\op{div}(Y)\right)-Y\left(\op{div}(X)\right),$$
and that verifies
$$\op{div}(fX)=f\op{div}(X)+X(f)$$
for any $X,Y\in\frak X(M)$ and $f\in C^\infty(M)$.
For details pertaining to divergence operators on an arbitrary
manifold, we refer the reader to \cite{GrabPoncinAuto}.

\begin{theo} Let Let $\zF:\cA_1\leftrightarrow\cA_2$ be an isomorphism of the Lie algebras $\cA_i$
of model reductive Atiyah algebroids $(A_i=\sZ_i\oplus\wA_i,[\![-,-]\!]_{\zw_i},\wzp_{i*}^0)$ with connected
reductive structure groups $G_i$ over connected manifolds $M_i$, $i\in\{1,2\}$, and let $\op{div}$ be a fixed
divergence on $M_1$. Then, there are
\begin{itemize}
\item a Lie algebroid isomorphism
$\zf^s:\wA_1\leftrightarrow\wA_2$ covering a diffeomorphism
$\zvf:M_1\leftrightarrow M_2$ and inducing a Lie algebra
isomorphism $\zF^s:\wcA_1\leftrightarrow\wcA_2$, \item a vector
bundle isomorphism $\zf^0:\sZ_1\leftrightarrow\sZ_2$ covering the
same diffeomorphism $\zvf$ and inducing a linear isomorphism
$\zF^0:\cZ_1\leftrightarrow\cZ_2$, \item a one-form $\zh$ on $M_1$
with values in $Z\frak g_1$ satisfying
$$d\zh=\zw_1-\zf^{0*}\zw_2\,,$$
\item an element $r\in Z\frak g_1$ representing a section of $\sZ_1$,
\end{itemize}
such that
\be\label{mf}
\zF(c+\wa)=\zF^0\left(c+\eta(\hat a)+\op{div}(\hat a)\cdot{r}\right)+\Phi^s(\wa)\,.
\ee
Conversely, every mapping of the form (\ref{mf}) with $\zF^s$, $\zF^0$, $\zh$, and $r$ satisfying the above conditions is a Lie algebra isomorphism.
\end{theo}

\begin{proof} We first show that
\be
\zF(c+\wa)=\zF^0(c+F(\wzp_{1*}\wa))+\zF^s(\wa),\label{SplitZF}\ee
where $\zF^0:\cZ_1\leftrightarrow \cZ_2$ and
$\zF^s:\wcA_1\leftrightarrow\wcA_2$ are the canonically
$\zF$-induced isomorphisms between the centers $\cZ_i$ and the Lie
algebras $\wcA_i$ of semisimple Atiyah algebroids, respectively,
see Proposition \ref{Phi0Prop}, Theorem \ref{TheoMain} and Theorem
\ref{TheoMain2}, and where $F:\frak X(M_1)\to \cZ_1$ is a linear
map. Indeed, let $\text{pr}_{\cZ_2}:\cZ_2\oplus\wcA_2\to\cZ_2$ be
the canonical projection. Since
$$\zF(c+\tilde a)=\left(\zF^0(c)+\text{pr}_{\cZ_2}(\zF(\tilde a))\right)+
\zF^s(\tilde a)\,,$$ it suffices to set
$F:=(\zF^0)^{-1}\text{pr}_{\cZ_2}\zF:\wcA_1\to\cZ_1$ and to prove
that this linear map factors through the quotient
$\wcA_1/\wcK_1\simeq{\frak X}(M_1)$, i.e. to show that
$\text{pr}_{\cZ_2}\zF$ vanishes on $\wcK_1\simeq \wcK_1^{(1)}$.
Indeed, when using Equation (\ref{LadBracketThroughIsom}), as well
as the fact that the $\sZ_i\oplus \tilde K_i$ are the kernel
{\small LAB}s of the Atiyah algebroids $\sZ_i\oplus \tilde A_i$,
we get
$$\zF\lp\wcK_1^{(1)}\rp=\zF\left((\cZ_1\oplus\wcK_1)^{(1)}\right)=(\cZ_2\oplus\wcK_2)^{(1)}=\wcK_2^{(1)}\,,$$
so that $\zF(\tilde{k})\in\wcA_2$, for each $\tilde{k}\in\wcK_1$.
It suffices now to put
$F(\wzp_{1*}\wa)=(\zF^0)^{-1}\left(\text{pr}_{\cZ_2}(\zF(\tilde
a))\right)$. Note that the map $\zF$ defined by (\ref{SplitZF}) is
always a linear isomorphism.

Now, we will show that it is a Lie algebra isomorphism if and only
if, for any vector field $X$, $F(X)=\eta(X)+\op{div}(X)\cdot r$,
with $\eta$ and $r$ as described in the theorem.

When applying Equations (\ref{LadBracketThroughIsom}) and (\ref{SplitZF}), as
well as the Lie algebra morphism property of $\zF^s$, we find that the Lie algebra morphism property
$$\zF\left([\![c+\wa,c'+\wa']\!]_{\zw_1}\right)=[\![\zF(c+\wa),\zF(c'+\wa')]\!]_{\zw_2}$$
reads
\bean\label{LaMorpCondZF}&\zF^0\lp
\hat a(c')-\hat a'(c)+\zw_1(\hat a,\hat a')+F(\wzp_{1*}[\wa,\wa']^{\,\tilde{}\;})\rp=\\
&(\wzp_{2*}\zF^s(\wa))\left(\zF^0(c'+F(\hat a'))\right)- (\wzp_{2*}\zF^s(\wa'))\left(\zF^0(c+F(\hat a))\right)
+\zw_2(\wzp_{2*}\zF^s(\wa),\wzp_{2*}\zF^s(\wa')).\nonumber\eean Note now that, since $\zF^s$ induces a Lie
algebra isomorphism $\tilde\zF^s=\zvf_*$ implemented by a diffeomorphism $\zvf$, we have
$$\wzp_{2*}\zF^s(\wa)=\zvf_*\wzp_{1*}\wa=\zvf_*\hat a\,,$$
and eventually combine the last upshot with Equation (\ref{IsomProp_ac}). This leads to
$$(\wzp_{2*}\zF^s(\wa))(\zF^0(c'+F(\hat a')))=(\zvf_*
\hat a)(\zF^0(c'+F(\hat a')))=\zF^0(\hat a(c'+F(\hat a'))),$$ so that the morphism condition
(\ref{LaMorpCondZF}) is equivalent to
\be\label{7}\omega_2(\zvf_*X,\zvf_*X')-\Phi^0\left(\omega_1(X,X')
-d F(X,X')\right)=0\,,\ee where
$$d F(X,X')=L_{X}(F(X'))-L_{X'}(F(X))
-F([X,X'])\,,$$ for all vector fields $X,X'$ of $M_1$. When decomposing $\omega_1=\sum_{\E}\omega_1^{\E}\frak
c_{1\E}$, $\omega_{2}=\sum_{\E}\omega_{2}^{\E}\frak c_{2\E}$, $F=\sum_{\E}F^{\E}\frak c_{1\E}$ in the
corresponding global frames, and when observing that
$$({\zvf}^*\zw_{2})(X,X')=\zw_{2}(\zvf_*X,\zvf_*X')\;\circ\zvf
\quad\mbox{and}\quad\zF^0\!\!\lp\sum_if^i\frak
c_{1i}\rp\;\circ\zvf=\sum_{ij}I^i_jf^j\frak c_{2i},$$ we can
rewrite (\ref{7}) in the form \be\label{7a}\xd
F_{i}^{m}=\omega_{1}^{m}-\sum_{\E}(I^{-1})_{\E}^m\zvf^*\omega_{2}^{\E}\,,\quad
m\in\{1,\dots,k\}\,. \ee Equation (\ref{7a}) implies that each
$F_{i}^m:\frak X(M_1)\to C^\infty(M_1)$ is a local, thus locally,
a differential operator. Indeed, if a vector field $Y$ vanishes in
a neighborhood of a point $m_1\in M_1$, it can be written in the
form $Y=\sum_k[X_k,X'_k]$ with vector fields $X_k,X'_k$ that
vanish in a neighborhood of $m_1$. It follows from (\ref{7a}) that
$F_{i}^m(Y)$ equals to
$$\sum_kF_{i}^m([X_k,X'_k])=\sum_k\left(X_k(F_{i}^m(X'_k))-X'_k(F_{i}^m(X_k))-
\lp\omega_{1}^m-\sum_{\E}(I^{-1})^m_{\E}\zvf^*\omega_{2}^{\E}\rp(X_k,X'_k)\right)\,,
$$ so that $F_{i}^m(Y)$ vanishes in a neighborhood of $m_1$.

Let $U_{\frak a}$, $\frak a\in\frak A$, be an open covering of $M_1$ by contractible charts. Since the 2-form $\zb^m=\omega_{1}^{m}-\sum_{\E}(I^{-1})_{\E}^m\zvf^*\omega_{2}^{\E}$ of
$M_1$ is closed, there is, for every $\frak a\in\frak A$, a one-form $\za^m_{\frak a}$ on $U_{\frak a}$ such
that $\zb^m|U_{\frak a}=\xd\za^m_{\frak a}$. The linear map $F^m|U_{\frak a}-\za^m_{\frak a}:\frak X(U_{\frak
a})\to C^\infty(U_{\frak a})$ is therefore a 1-cocycle of the Chevalley-Eilenberg cohomology of vector fields
represented upon functions, i.e., see \cite{GrabPoncinAuto}, \cite{NP}, $F^m|U_{\frak a}-\za^m_{\frak
a}=r^m_{\frak a}\op{div}+\gamma^m_{\frak a}$, for some $r^m_{\frak a}\in\R$ and some closed 1-form
$\gamma^m_{\frak a}\in\zW^1(U_{\frak a})$. In other words, $F^m|U_{\frak a}=r^m_{\frak
a}\op{div}+\eta^m_{\frak a}$, where the one-form $\eta^m_{\frak a}=\gamma^m_{\frak a}+\za^m_{\frak a}$. Since
on an intersection $U_{\frak a}\cap U_{\frak a'}$ the constants $r^m_{\frak a}$ and $r^m_{\frak a'}$ must
coincide with a single $r^m$ (as $M_1$ is connected), the one-forms $\eta^m_{\frak a}$ and $\eta^m_{\frak a'}$
coincide as well. Hence $F^m=r^m\op{div}+\eta^m$. Since $\op{div}$ is a cocycle, we have
$d\eta^m=\zb^m$.\end{proof}

It is clear that if we choose $r=0$ we get an isomorphism of Lie algebroids $\zF_0:\cA_1\ni
c+\tilde{a}\leftrightarrow \zF^0(c+\zh(\hat a))+\zF^s(\tilde a)\in \cA_2$. It follows that $\zF=\zF_0\zF_1$,
where $\zF_1:\cA_1\ni c+\tilde{a}\leftrightarrow c+\tilde{a}+\op{div}(\hat a)\cdot r\in\cA_1$ is an
automorphism of $\cA_1$. When identifying the Atiyah algebroids with the model algebroids, we get Theorem
\ref{t0}.

\vskip1cm \noindent Janusz GRABOWSKI\\ Polish Academy of Sciences, Institute of Mathematics\\ \'Sniadeckich 8,
P.O. Box 21, 00-956 Warsaw,
Poland\\Email: jagrab@impan.pl \\

\noindent Alexei KOTOV\\
University of Luxembourg, Mathematics Research Unit\\ 162A, avenue de la Fa\"iencerie, L-1511 Luxembourg City,
Grand-Duchy of Luxembourg\\Email: alexei.kotov@uni.lu\\

\noindent Norbert PONCIN\\
University of Luxembourg, Mathematics Research Unit\\ 162A, avenue
de la Fa\"iencerie, L-1511 Luxembourg City, Grand-Duchy of
Luxembourg\\ E-mail: norbert.poncin@uni.lu

\end{document}